\documentclass{amsart}[11pt]
\usepackage{mathrsfs, amsmath,amssymb}
\usepackage[utf8x]{inputenc}
\usepackage{fancyhdr,amsfonts,stmaryrd}
\usepackage{textcomp}
\newtheorem{theo}{Theorem}[section]

\newtheorem{cor}[theo]{Corollary}
\newtheorem{lem}[theo]{Lemma}
\newtheorem{pro}[theo]{Proposition}
\newtheorem{Def}[theo]{Definition}
\newtheorem{rem}[theo]{Remark}
\newtheorem{Not}[theo]{Notations}
\newtheorem{cons}[theo]{Consequence}
\newtheorem{Ex}[theo]{Example}
\def\pf{{ \textbf{\textit{Proof.}}}\quad}

\newcommand{\N}{\mathbb{N}}
\newcommand{\Q}{\mathbb{Q}}
\newcommand{\R}{\mathbb{R}}

\newcommand{\E}{\mathbb{E}}
\newcommand{\p}{\mathbb{P}}

\newcommand{\boF}{\mathcal{F}}

\newcommand{\boI}{\mathcal{I}}
\newcommand{\boL}{\mathcal{L}}

\newcommand{\boU}{\mathcal{U}}
\newcommand{\boH}{\mathcal{H}}
\newcommand{\boS}{\mathcal{S}}

\title{Stochastic Analysis for obtuse random walks}

\author[U. Franz]{Uwe Franz}
\address{UFC, Département de mathématiques de Besançon, 16 route de Gray, 25 030 Besançon cedex, France}
\email{uwe.franz@univ-fcomte.fr}

\author[T. Hamdi]{Tarek Hamdi}
\address{IPEST, Universit\'e de Carthage\\ 2078 La Marsa \\ Tunisie}
\email{tarek.hamdi@ipest.rnu.tn}

\begin{document}

\begin{abstract}
We present a construction of the basic operators of stochastic
analysis (gradient and divergence) for a class of discrete-time normal martingales called obtuse random walks. 
The approach is based on the chaos representation property and discrete multiple
stochastic integrals. We show that these operators 
satisfy similar identities as in the case of the Bernoulli randoms walks. We prove a Clark-Ocone-type predictable representation 
formula, obtain two covariance identities and derive a deviation inequality.
We close the exposition by an application to option hedging in discrete time. 
\end{abstract}
\maketitle

\section{Introduction}

A celebrated Theorem of Wiener \cite{W} (who introduced the terms 'homogeneous chaos' and 'polynomial chaos' in that paper) asserts that the chaotic representation property (hereafter CRP) holds
for the Brownian motion. This property says that any square integrable random variable measurable with respect to a Brownian motion $X$ can be expressed as an 
orthogonal sum of multiple stochastic integrals with respect to $X$. 
The main feature of this property is that it gives rise to an isometry between the Fock space and the $L^2$-space associated with 
this Brownian motion. In particular, the $n$-th Wiener chaos is identified to an element of the Fock space for any $n \geq 1$, opening the way to set up anticipating and non commutative stochastic calculus.      
Actually, Wiener chaoses were slightly different from the modern ones, introduced later by Itô \cite{I} and studied by Meyer \cite{M} for an interesting class of 
martingales called normal martingales (including Brownian motion and the compensated Poisson process). In addition to the martingale property, these processes, say $X = (X^1,\dots,X^d)$, are specified by the requirement  
\begin{equation*}
 <X^i,X^j>_t=\delta^{ij}t.
\end{equation*}
Besides, Meyer noticed that any normal martingale gives rise to an isometry between the Fock space and its $L^2$-space. When this isometry is an isomorphism 
of Hilbert spaces, then it leads to structure equations (hereafter SE) i.e. we have
\begin{equation*}
 [X^i,X^j]=\delta^{ij}t+\int_0^t\sum_k(\Phi_k^{ij})_sdX_s
\end{equation*}
for some predictable processes $(\Phi_k^{ij})$, and furthermore, $X$ enjoys the PRP. However, this is far from being a sufficient condition unlike its analog in the discrete-time
setting introduced and developed by Attal and Emery in \cite{AE}. In fact, in the discrete time case, Attal and Emery showed that
SEs are necessary and sufficient for the CRP to hold (and also for the predictable 
representation property (PRP)). We are thus led to the so-called obtuse random walks which are a class of $d$-dimensional normal martingales such that the sequence of their
increments $(\Delta X_n)_{n\geq0}$ take $d+1$ values for each $n$. This fact translates that the filtration $(\boF_n)_{n\geq0}$ generated by $X$ is of
multiplicity $d+1$ (i.e. we move from $\boF_n$ to $\boF_{n+1}$ by decomposing of each atom of
$\boF_n$ to $d+1$ atoms) c.f. \cite{AE}. For $d=1$, these random walks reduce to the Bernoulli process  which was used in \cite{P} to deal with the CRP 
and to define discrete multiple stochastic integrals with respect to a discrete-time normal martingale with i.i.d. sequences of increments.
 
In this paper, we shall focus on discrete time normal $d$-dimensional martingales (see definition below). 
Our main concern is generalizing the stochastic analysis for Bernoulli random walks (see \cite{P}) to the obtuse random walks ($d>1$) with 
not necessarily independent increments. First, we present a construction of the stochastic integral of predictable square-integrable processes and the 
associated multiple stochastic integrals of symmetric functions on $\N^n\ (n ≥ 1)$, with respect to such martingales.
Indeed, these iterated stochastic integrals give an isometry between
this $L^2$-space and the Fock space. 
Next, we present a construction of the basic operators of stochastic analysis (gradient and divergence, \cite{P}). We give a probabilistic interpretation of the gradient operator and prove that the divergence operator coincides with the stochastic integral on square summable predictable processes.
These operators are used to derive a Clark-Ocone-type predictable representation formula and also to prove a deviation inequality. 
Finally, we apply the tools developed in this paper to discrete market models in order to obtain explicit expressions for hedging strategies \cite{KO,P}.

One motivation for this paper is of course to get a better intuition for the much more difficult continuous-time case, like in \'Emery's paper \cite{E}. 
But, so far we have no new results in this direction. Another reason for our study is the hope to apply our tools of multi-dimensional discrete-time stochastic 
analysis to non-commutative discrete-time stochastic calculus, e.g. in models of repeated quantum interaction \cite{dhahri+attal}.

This paper is organized as follows. In section 2 we present a construction of the stochastic integral of predictable square-integrable processes 
with respect to a normal martingale. 
In the next section we construct the associated multiple stochastic integrals of functions $f_n$ that are symmetric in $n$ variables.
In section 4 we present a characterization of obtuse random walks in discrete-time setting.
The proof of the CRP is reviewed in section 5. A gradient operator $D$ acting by finite differences is introduced in section 6 in connection with multiple stochastic
integrals and is used in section 7 to state a Clark-Ocone-type predictable representation formula. The divergence operator $\delta$ is defined in section 8 as an
extension of the discrete-time stochastic integral and we shall prove that it is the adjoint of $D$. The Ornstein-Uhlenbeck semi-group is used in section 9 to
express a covariance identity. In section 10 we prove a deviation inequality for functionals of obtuse random walks.
The last section is devoted to present a complete market model in discrete time as an application of the Clark-Ocone Formula.

\section{Discrete stochastic integrals}

Consider a discrete $d$-dimensional process $Y=(Y^1,...,Y^d)$ on a probability space $(\Omega,\boF,\p).$ Let $(\boF_n)_{n\geq0}$ denote the filtration
generated by $(Y_n)_{n\in\N}$ and $\boF_{-1}=\{\varnothing,\Omega\}$.
Recall that a $d$-dimensional integrable process is said to be an $\boF_n$-martingale if each coordinate is so.

We recall also (see for example \cite{AE}) that a d-dimensional martingale $(Y_0+\ldots+Y_n)_{n\in\N}$ is said to be normal martingale if for any $n\geq0$, 
\begin{equation*}
 \E[Y_n^iY_n^j\,|\boF_{n-1}]=\delta^{ij},\quad i,j\in\{1,...,d\},
\end{equation*}
which can be written as
\begin{equation*}
 \E[Y_n\otimes Y_n\,|\boF_{n-1}]=I_n,
\end{equation*}
where $\otimes$ is the Kronecker tensor product of the vector $Y_n$ by itself.

In the sequel, we denote $<x,y>$ the inner product of $x$ and $y$ in $\R^d$, and we make following assumptions on $Y=(Y^1,...,Y^d)$:
\begin{equation*}
 \E[Y_n^i\,|\boF_{n-1}]=0\quad and\quad \E[Y_n^iY_n^j\,|\boF_{n-1}]=\delta^{ij}.
\end{equation*}
These assumptions imply that the process $(Y_0 +...+ Y_n)_{n\geq0}$ is a normal martingale in the discrete time.
\begin{Def}
 Let $U=(U^1,...,U^d)$ be a uniformly bounded sequence of random variables with finite support in $\N$ (i.e. $U_n=0_{\R^d}$ except perhaps a finite number of indices).
 The stochastic integral of $U$ with respect to $Z$ is defined as
\begin{equation*}
\mathcal{I}(U)=\sum_{k=1}^{d}\sum_{n=0}^{\infty}U_n^kY_n^k=\sum_{n=0}^{\infty}<U_n,Y_n>.
\end{equation*}

\end{Def}
Now we recall that:
\begin{Def}
 A stochastic process $(X_n)_n$ is said to be predictable process with respect to $(\boF_n)_n$ if $X_n$ is $\boF_{n-1}$-measurable for each $n$.
\end{Def}
Then one has the following result.
\begin{pro}
The stochastic integral extends to square-integrable predictable processes via the (conditional) isometry formula:
\begin{equation*}
 \forall n\in \N,\quad \E[\lvert\boI( 1_{[n,\infty)}U)\lvert^2\,|\boF_{n-1}]=\E[\Arrowvert 1_{[n,\infty)}U\Arrowvert^2\,|\boF_{n-1}]
\end{equation*}
where $1_{[n,\infty)}U$ denotes  the process $(0,...,0,U_n,U_{n+1},...).$
\end{pro}
 \pf
Let $U,V$ be bounded predictable processes with finite support in $\N,$ we have
 \begin{eqnarray*}
  \E\left[\sum_{k=n}^{\infty}<U_k,Y_k>\sum_{l=n}^{\infty}<V_l,Y_l>\,|\boF_{n-1}\right]
&=&\E\left[\sum_{i,j=1}^{d}\sum_{k,l=n}^{\infty}U_k^iY_k^iV_l^jY_l^j\,|\boF_{n-1}\right]\\
&=&\sum_{i,j=1}^{d}\sum_{k=n}^{\infty}\E\left[\E[U_k^iY_k^iV_k^jY_k^j\,|\boF_{k-1}]\,|\boF_{n-1}\right]\\
&&+\sum_{i,j=1}^{d}\sum_{n\leq k<l}\E\left[\E[U_k^iY_k^iV_l^jY_l^j\,|\boF_{l-1}]\,|\boF_{n-1}\right]\\
&&+\sum_{i,j=1}^{d}\sum_{n\leq l<k}\E\left[\E[U_k^iY_k^iV_l^jY_l^j\,|\boF_{k-1}]\,|\boF_{n-1}\right]\\
&=&\sum_{i,j=1}^{d}\sum_{k=n}^{\infty}\E\left[U_k^iV_k^j\E[Y_k^iY_k^j\,|\boF_{k-1}]\,|\boF_{n-1}\right]\\
&&+2\sum_{i,j=1}^{d}\sum_{n\leq l<k}\E\left[U_k^iY_k^iV_l^j\E[Y_l^j\,|\boF_{k-1}]\,|\boF_{n-1}\right]\\
&=&\sum_{i=1}^{d}\sum_{k=n}^{\infty}\E\left[U_k^iV_k^j\,|\boF_{n-1}\right]\\
&=&\E\left[\sum_{k=n}^{\infty}<U_k,V_k>\,|\mathcal{F}_{n-1}\right].
 \end{eqnarray*}
$\hfill\square$

\section{Discrete multiple stochastic integrals }

This section is devoted to the construction and to state the main properties of the multiple stochastic
integrals of symmetric functions on $\N^r, r \geq 1.$
We denote 
\begin{equation*}
 \Delta_r=\{(i_1,...,i_r)\in\N^r,i_l\neq i_k, 1\leq l<k\leq r\}
\end{equation*}
 and
\begin{equation*}
 \begin{array}{lccc}
  f_r:&\Delta_r&\longrightarrow&\R^{d^r}
\\&(i_1,...,i_r)&\longmapsto&\left( f_r^{k_1,...,k_r}(i_1,...,i_r)\right)_{1\leq k_1,...,k_r\leq d} 
\end{array}
\end{equation*}

a symmetric function in $r$ variables. 

Given $f_1\in l^2(\N)$ we let
\begin{equation*}
 \boI^1(f_1)=\boI(f_1)=\sum_{n=0}^{\infty}<f_1(n),Y_n>.
\end{equation*}

\begin{Def}
 For $\ r\geq1$, the multiple stochastic integral of $f_r\in L^2(\Delta_r,\R^{d^r})$ with respect to the normal martingale $(Y_0+\ldots+Y_n)_{n\geq0}$ is defined by 
\begin{equation*}
 \boI^r(f_r)=\sum_{k_1,...,k_r=1}^d\ \sum_{(i_1,...,i_r)\in\Delta_r}f_r^{k_1,...,k_r}(i_1,...,i_r)Y_{i_1}^{k_1}...Y_{i_r}^{k_r}.
\end{equation*}
\end{Def}
\begin{rem}
 We take $\Delta_0=\{0\}$, $ L^2(\Delta_0,\R)=\R$ and define $\boI^0(f_0)\equiv f_0,\ f_0\in\R.$ 
\end{rem}
The following result gives a recurrence relation for multiple stochastic integral.
\begin{pro}
Let $r\geq 1$, we have
\begin{equation*}
 \boI^r(f_r)=r\sum_{k_r=1}^d\ \sum_{i_r=0}^{\infty}\boI^{r-1}\left( f_r^{k_r}(*,i_r)1_{{\text{\textlbrackdbl}0,i_r-1\text{\textrbrackdbl}}^{r-1}}(*)\right) Y_{i_r}^{k_r},
\end{equation*}
where 
\begin{equation*}
 \begin{array}{lccc}
f_r^{k}(*,i):&\Delta_{r-1}&\longrightarrow&\R^{d^{r-1}}
\\&(i_1,...,i_{r-1})&\longmapsto&\left( f_r^{k_1,...,k_{r-1},k}(i_1,...,i_{r-1},i)\right)_{1\leq k_1,...,k_{r-1}\leq d} 
\end{array}
\end{equation*}

\end{pro}
\pf We write
\begin{eqnarray*}
 \boI^r(f_r)&=&r!\sum_{k_1,...,k_r=1}^d\ \sum_{i_r=0}^{\infty}\sum_{0\leq i_{r-1}\leq i_r}...\sum_{0\leq i_{1}\leq i_2}f_r^{k_1,...,k_r}(i_1,...,i_r)Y_{i_1}^{k_1}...Y_{i_r}^{k_r}\\
&=&r\sum_{k_r=1}^d\ \sum_{i_r=0}^{\infty}\left( (r-1)!\sum_{k_1,...,k_{r-1}=1}^d\ \sum_{0\leq i_{r-1}\leq i_r}...\right. 
\\ &&...\left. \sum_{0\leq i_{1}\leq i_2}f_r^{k_1,...,k_r}(i_1,...,i_r)Y_{i_1}^{k_1}...Y_{i_{r-1}}^{k_{r-1}}\right)Y_{i_r}^{k_r}.
\end{eqnarray*}
$\hfill\square$

The next proposition states an isometry formula
\begin{pro}
Let $r,s\geq1,$ and consider 
\begin{equation*}
 f_r=(f_r^{k_1,...,k_r}(i_1,...,i_r))_{1\leq k_1,...,k_r\leq d}\in L^2(\Delta_r,\R^{d^r}), 
\end{equation*}
\begin{equation*}
g_s=(g_s^{t_1,...,t_s}(j_1,...,j_s))_{1\leq t_1,...,t_s\leq d}\in L^2(\Delta_s,\R^{d^s}).
\end{equation*}
 We have
\begin{equation*}
 \E\left[\boI^r(f_r)\boI^s(g_s) \right]= 1_{\{s=r\}}r!\sum_{k_1,...,k_r=1}^d<f_r^{k_1,...,k_s},g_s^{k_1,...,k_s}>
\end{equation*}
\end{pro}
\pf 
\begin{eqnarray*}
 &&\E\left[\boI^r(f_r)\boI^s(g_s) \right]
\\&&= \displaystyle\sum_{
\begin{array}{l}
\scriptstyle k_1,...,k_r=1\\
\scriptstyle t_1,...,t_s=1
                    \end{array}
}^d \displaystyle\sum_{\begin{array}{l}
\scriptstyle(i_1,...,i_r)\in\Delta_r\\
\scriptstyle(j_1,...,j_s)\in\Delta_s
                    \end{array}
}
f_r^{k_1,...,k_r}(i_1,...,i_r)g_s^{t_1,...,t_s}(j_1,...,j_s)\E[Y_{i_1}^{k_1}...Y_{i_r}^{k_r}Y_{j_1}^{t_1}...Y_{j_s}^{t_s}]\\
&&=(r!)^2\displaystyle\sum_{
\begin{array}{l}
\scriptstyle k_1,...,k_r=1\\
\scriptstyle t_1,...,t_s=1
                    \end{array} 
}^d
\sum_{\begin{array}{l}
\scriptstyle0\leq i_1<...<i_r\\
\scriptstyle0\leq j_1<...<j_s
                    \end{array}
}
f_r^{k_1,...,k_r}(i_1,...,i_r)g_s^{t_1,...,t_s}(j_1,...,j_s)\E[Y_{i_1}^{k_1}...Y_{i_r}^{k_r}Y_{j_1}^{t_1}...Y_{j_s}^{t_s}]
\end{eqnarray*}
Note that if $r=s$ and $0\leq i_1<...<i_r$ and $0\leq j_1<...<j_r$ we have 
\begin{equation*}
 \E[Y_{i_1}^{k_1}...Y_{i_r}^{k_r}Y_{j_1}^{t_1}...Y_{j_s}^{t_s}]=1_{\{i_1=j_1,...i_r=j_r\}}1_{\{k_1=t_1,...k_r=t_r\}},
\end{equation*}
hence we get
\begin{eqnarray*}
 &&\E\left[\boI^r(f_r)\boI^r(g_r) \right]
\\&&=r!\sum_{
\begin{array}{l}
\scriptstyle k_1,...,k_r=1\\
\scriptstyle t_1,...,t_s=1
                    \end{array}
}^d\ \sum_{(i_1,...,i_r)\in\Delta_r}f_r^{k_1,...,k_r}(i_1,...,i_r)g_r^{k_1,...,k_r}(i_1,...,i_r)1_{\{k_1=t_1,...k_r=t_r\}}.
\end{eqnarray*}
If $r<s$ then there necessarily exists $k\in\{1,...,s\}$ such that $j_k\notin\{i_1,...,i_r\}$ thus 
\begin{equation*}
 \E[Y_{i_1}^{k_1}...Y_{i_r}^{k_r}Y_{j_1}^{t_1}...Y_{j_s}^{t_s}]=0.
\end{equation*}

$\hfill\square$

\section{Obtuse random walks}

Let us recall briefly the canonical construction of discrete-time normal martingales with values in $\R^d$.
Consider a normal martingale $(Y_0 +\ldots+ Y_n)_{n\geq0}$ such that, for each $n$, $Y_n$ takes $d+1$ values 
$v_0(n),...,v_d(n)$ conditionally to $\boF_{n-1}$. Let $\p$ be any probability measure on the set $\Omega=\{0,...,d\}^{\N}$ that assigns
strictly positive probability $p^i_n$ to each $v_i(n)$ where $(v_i(n))_n$ and $(p_n^i)_n$ are predictable processes.
$(\boF_n)_{n\geq0}$ denote the filtration generated by $(Y_n)_{n\in\N}$ i.e.
\begin{equation*}
 \boF_n=\sigma(Y_0,...,Y_n),\quad n\in\N.
\end{equation*}
 We introduce the coordinate maps 
\begin{equation*}
 \begin{array}{lllc}
   X_n:&\Omega&\longrightarrow&\{0,1,...,d\}\\
&w&\longmapsto&w_n
  \end{array}
\end{equation*}

For 
\begin{equation*}
 w=(w_0,w_1,...,w_n,.....)\in\Omega,
\end{equation*}
we write $Y_n(w)=v_{X_n(w)}$ which yields
\begin{equation*}
 \boF_n=\sigma(X_0,...,X_n),\quad n\in\N.
\end{equation*}
Hence we have 
\begin{equation*}
 p_n^i=\p(Y_n=v_i(n)\,|\boF_{n-1})=\p(X_n=i\,|\boF_{n-1}),\quad n\in\N.
\end{equation*}
Let 
\begin{equation*}
 c_i^j(n)=p^i_nv_i^j(n),\quad n\in\N,\,i\in\{0,...,d\}\  \text{and} \ j\in\{1,...,d\}.
\end{equation*}
\begin{pro}\label{aa}
$\forall n\in\N,\,\forall j,l\in\{1,...,d\}$, we have 
\begin{equation*}
 \sum_{i=0}^dc_i^j(n)=0,
\end{equation*}
and
$$\sum_{i=0}^dc_i^j(n)v_i^l(n)=\delta^{jl}.$$
\end{pro}
\pf We write
\begin{equation*}
 \sum_{i=0}^dc_i^j(n)=\E[Y_n^j\,|\boF_{n-1}]=0,
\end{equation*}
and 
\begin{equation*}
 \sum_{i=0}^dc_i^j(n)v_i^l(n)=\E[Y_n^jY_n^l\,|\boF_{n-1}]=\delta^{jl}.
\end{equation*}
$\hfill\square$

Recall that the filtration $(\boF_n)_{n\geq0}$ is said to be of multiplicity $d+1$ if each $\boF_n$ is finite and
each atom of $\boF_n$ contains exactly $d+1$ atoms of $\boF_{n+1}$.
The following result gives a characterization of normal martingales which satisfy the CRP (c.f. \cite{AE} for a proof and for more details).
\begin{theo}\label{p}
 Let $(Y_0 +\ldots+ Y_n)_{n\geq0}$ be a d-dimensional normal martingale, the following assertions are equivalent
\begin{enumerate}
 \item The filtration multiplicity is bounded from above by $d+1$.
\item The filtration multiplicity is exactly $d+1$.
\item $(Y_0 +\ldots+ Y_n)_{n\geq0}$ satisfies a SE 
\begin{equation*}
 Y_n^iY_n^j=\delta^{ij}+\sum_{k=1}^d\Phi_{ij}^k(n)Y_n^k,
\end{equation*}
where $\Phi_{ij}^k$ are $d^3$ predictable processes.
\item $(Y_0 +\ldots+ Y_n)_{n\geq0}$ has the PRP.
\item $(Y_0 +\ldots+ Y_n)_{n\geq0}$ has the CRP. 
\end{enumerate}
\end{theo}
\begin{Def}
 An obtuse random walk is a process that satisfies the equivalent condition of Theorem \ref{p}.
\end{Def}

Note that the values $(v_i(n))_{0\leq i\leq d}$ of \,$Y_n$ and their probabilities $(p_n^i)_{0\leq i\leq d}$ are related to the coefficients of the SE by
\begin{equation*}
 \Phi(n)=\sum_{i=0}^dp_n^iv_i^*(n)\otimes v_i(n)\otimes v_i(n).
\end{equation*}

\section{Chaotic representation property}

Assume now that the filtration $(\boF_n)_n$ generated by $(Y_n)_{n\in\N}$ has a multiplicity equal to $d+1$.  
Let $L^0(\Omega,\boF_n)$ the space of $\boF_n$-measurable random variables, it has finite dimension equal to $(d+1)^{n+1}.$\\
For $N\in\N,$ we denote
\begin{equation*}
 \boI_N^r(f_r)=\boI^r\left( f_r1_{{\text{\textlbrackdbl}0,N\text{\textrbrackdbl}}^{r}}\right). 
\end{equation*}
Note that if $r>N+1$, then $\boI_N^r(f_r)=0.$
\begin{pro}
For all $r\geq1,$
 \begin{equation*}
  \boI_N^r(f_r)=\E\left[ \boI^r(f_r)\,|\boF_N\right].
 \end{equation*}
\end{pro}
\pf Let $0\leq i_1<...<i_r\ \in\Delta_r,$ if $i_r>N$ we have
\begin{equation*}
 \E[Y_{i_1}^{k_1}...Y_{i_r}^{k_r}]=\E[\E[Y_{i_1}^{k_1}...Y_{i_r}^{k_r}]\,|\boF_{r-1}]=0.
\end{equation*}
As a result 
\begin{equation*}
 \E[\boI^r(f_r)]=0,\ \ \forall r\geq1
\end{equation*}
and the process $(\boI_k^r(f_r))_{k\in\N}$ is a discrete-time martingale. 
$\hfill\square$

\begin{cor} For
 $ 0\leq N\leq r$,\\ 
$\boI^r(f_r)$ is $\boF_N-$measurable if and only if $f_r1_{{\text{\textlbrackdbl}0,N\text{\textrbrackdbl}}^{r}}=f_r.$
\end{cor}
\pf The sufficiency is obvious. The necessity is a consequence of
\begin{equation*}
 \boI_N^r(f_r)=\E\left[ \boI^r(f_r)\,|\boF_N\right]=\boI^r(f_r),
\end{equation*}
and of the isometry formula.
$\hfill\square$

\begin{Def}
 Let $\boH_0=\R$ and for $n\geq1$, we denote $\boH_n$ the subspace of $L^2(\Omega)$ made of stochastic integrals of order $n\geq1$
\begin{equation*}
 \boH_n=\{\boI^n(f_n), f_n\in L^2(\Delta_n,\R^{d^n})\}.
\end{equation*}
\end{Def}
\begin{pro}\label{crp}
 $\forall n\in\N,$
\begin{equation*}
 L^0(\Omega,\boF_n)\subset \boH_0\oplus...\oplus \boH_{n+1}.
\end{equation*}
\end{pro}
\pf For $0\leq r\leq n+1$, we have dim$\,L^0(\Omega,\boF_n)\cap\boH_r=\binom{n+1}{r}d^r $. More precisely,
\begin{equation*}
 \{Y_{i_1}^{k_1}...Y_{i_r}^{k_r}:\ 0\leq i_1<...<i_r\leq n,\ 1\leq k_1,...,k_r\leq d\}
\end{equation*}
form an orthonormal basis. By orthogonality of the subspaces $\boH_r$ we have
\begin{equation*}
 L^0(\Omega,\boF_n) = (\boH_0\oplus...\oplus \boH_{n+1})\cap L^0(\Omega,\boF_n).
\end{equation*}

$\hfill\square$

\begin{cons}
 Any element $F\in L^2(\Omega,\boF_n)$ can be written as
\begin{equation*}
 F=\E[F]+\sum_{r=1}^{n+1}\boI_n^r(f_r).
\end{equation*}
\end{cons}
\begin{Def}
 We denote $\boS$ the linear space spanned by multiple stochastic integrals 
\begin{equation*}
 \boS=\left\lbrace \bigcup_{n=0}^{\infty}\boH_n\right\rbrace =\left\lbrace \sum_{r=0}^n\boI^r(f_r), with\  f_r\in L^2(\Delta_r,\R^{d^r})\ symmetric\right\rbrace.
\end{equation*}
The completion of $\boS$ in $L^2(\Omega)$ is denoted by the direct sum 
\begin{equation*}
 \bigoplus_{n=0}^{\infty}\boH_n.
\end{equation*}
\end{Def}
The next result establishes the CRP for normal martingales under the assumption that we move from $\boF_n$ to $\boF_{n+1}$ by decomposition of each atom of
$\boF_n$ to $d+1$ atoms (with measure $>0$), see \cite{AE,M}. 
\begin{theo}
\begin{equation*}
 L^2(\Omega)=\bigoplus_{n=0}^{\infty}\boH_n.
\end{equation*}
\end{theo}
\pf It suffices to show that $\boS$ is dense in $L^2(\Omega)$. To this end, let $F$ be a bounded random variable, then 
Proposition \ref{crp} shows that, $\E[F\,|\boF_n]\in\boS$.
But, the martingale $(\E[F\,|\boF_n])_{n\in\N}$ converges a.s. and in $L^2(\Omega)$ to $F$, and we are done. 
$\hfill\square$

\section{Gradient operator}

\begin{Def}
 The gradient operator $D:\boS\longrightarrow L^2(\Omega\times\N,\R^d)$ is defined by
\begin{equation*}
 D_k^j(\boI^r(f_r))=r\boI^{r-1}\left( f_r^{j}(*,k)1_{\Delta_r}(*,k)\right),\quad\ k\in \N\ and\ j\in\{1,...,d\}.
\end{equation*}
\end{Def}
\begin{pro}
 The gradient operator is continuous on the chaos $\boH_r$.
\end{pro}
\pf We have
\begin{eqnarray*}
 \|D_k\boI^r(f_r)\|_{L^2(\Omega,\R^d)}^2&=&\sum_{j=1}^d\|D_k^j\boI^r(f_r)\|_{L^2(\Omega,\R)}^2\\
&=&\sum_{j=1}^dr^2\|\boI^{r-1}\left( f_r^{j}(*,k)1_{\Delta_r}(*,k)\right)\|^2\\
&=&\sum_{j=1}^dr^2(r-1)!\|f_r^{j}(*,k)\|^2_{L^2(\Delta^{r-1})}\\
&=&rr!\|f_r(*,k)\|^2_{L^2(\Delta^{r-1})}.
\end{eqnarray*}
$\hfill\square$

\begin{pro}\label{P}
 Let $F\in\boS$ be $\boF_n$-measurable, then for any $ k>n,$ one has
\begin{equation*}
 D_k^jF=0,\quad\ j\in\{1,...,d\}.
\end{equation*}
\end{pro}
\pf
We write
\begin{equation*}
 F=\E[F]+\sum_{r=1}^{n+1}\boI_n^r(f_r)
\end{equation*}
Then for $k>n$ we have
\begin{equation*}
 D_k^j\left(\boI_n^r(f_r) \right)=r\boI^{r-1}(f_r^j(*,k)1_{[0,n]^r}^j(*,k)1_{\Delta_r}(*,k))=0,\quad \forall j\in\{1,...,d\}.
\end{equation*}

$\hfill\square$

\begin{rem}
 By the Clark-Ocone formula derived in the next section, the converse of this proposition is also true i.e. if $F\in\boS$ is such that
\begin{equation*}
 D_k^jF=0,\ \forall k>n\ \text{and}\ \forall j\in\{1,...,d\},
\end{equation*}
then $F$ is $\boF_n$-measurable.
\end{rem}
\begin{Not}
 Let $\widetilde{f}_n\in L^2(\R^{n},\R^{d^n})$ denote the symmetrization of $f_n\in L^2(\Delta_{n},\R^{d^n}),$ given by
\begin{equation*}
 \widetilde{f}_n^{i_1,...,i_n}(t_1,...,t_n)=\frac{1}{n!}\sum_{\sigma\in\boS_n}f_n^{i_{\sigma(1)},...,i_{\sigma(n)}}(t_{\sigma(1)},...,t_{\sigma(n)}),
\quad 1\leq i_1,...,i_n\leq d.
\end{equation*}
In particular, for $(s_1,...,s_r)\in\Delta_r $, we have
\begin{equation*}
 \widetilde{1}_{\{s_1,...,s_r\}}^{i_1,...,i_r}(t_1,...,t_r)=\frac{1}{r!}1_{\{\{s_1,...,s_r\}=\{t_1,...,t_r\}\}}e_{i_1}\otimes...\otimes e_{i_r},
\end{equation*}
where $(e_1,...,e_d)$ denotes the canonical basis of $\R^d.$ 
\end{Not}
\begin{pro}\label{S} For any $r\geq1$, we have
 \begin{equation*}
 \boI^r\left(\widetilde{1}_{(s_1,...,s_r)}^{i_1,...,i_r} \right)=Y^{i_1}_{s_1}...\,Y^{i_r}_{s_r}. 
 \end{equation*}
\end{pro}
As a result an orthonormal basis of $L^2(\Omega,\boF_n)$ is given by
\begin{equation*}
 \left\lbrace Y^{i_1}_{s_1}...\,Y^{i_r}_{s_r}:\ 0\leq s_1<...<s_r\leq n,\ 1\leq i_1,...,i_r\leq d \right\rbrace.
\end{equation*}
Hence we can write
\begin{pro}\label{s} For any $r\geq1$,
 \begin{equation*}
  D_k^j\left(Y^{i_1}_{s_1}...Y^{i_r}_{s_r} \right) =
\begin{cases}
\begin{array}{ll}

\delta^{ji_t} Y^{i_1}_{s_1}...\check{Y}_{s_t}^{i_t}...Y^{i_r}_{s_r}& if\ k=s_t,\ \ t\in\{1,...,r\}\\
0& if\ k\notin (s_1,...s_r)
   \end{array}
\end{cases}.
 \end{equation*}
where $\check{Y}_{s_t}^{i_t} $ denotes that the factor $Y_{s_t}^{i_t}$ should be omitted in the product.
\end{pro}
\pf
Using Proposition \ref{S}, one can see that
\begin{eqnarray*}
 D_k^j\left(Y^{i_1}_{s_1}...Y^{i_r}_{s_r} \right) &=&D_k^j\left( \boI^r\left(\widetilde{1}_{(s_1,...,s_r)}^{i_1,...,i_r}\right)\right)\\
  &=&r\boI^{r-1}\left(\frac{1}{r!}1_{\{j\in(i_1,...,i_r)\}}1_{\{(s_1,...,s_r)=(*,k)\}}1_{\Delta_r}(*,k)e_{i_1}\otimes...\otimes e_{i_{r-1}}\otimes e_j\right) \\
\end{eqnarray*}
$\hfill\square$\\
The following result gives the probabilistic interpretation of $D_k^j$ as a finite difference operator in the case of 
discrete time random walks with i.i.d. increments. 
\begin{pro}
 For any $F\in\boS$, one has
\begin{equation*}
 D_k^jF(w)=\sum_{i=0}^dc_i^j(k)F(w_i^k)
\end{equation*}
where we recall that (c.f. Section 4) 
$$c_i^j(k)=\p(X_k=i\,|\boF_{k-1})=Y_k^j(w_i^k)$$ and $$w_i^k=(w_1,...,w_{k-1},i,w_{k+1},...).$$
\end{pro}
\pf It suffices to consider $F=Y_{s_1}^{i_1}...Y_{s_r}^{i_r}.$ By Proposition \ref{s} 
\begin{equation*}
 D_k^jF=\begin{cases}
\begin{array}{ll}

\delta^{ji_t} Y^{i_1}_{s_1}...\check{Y}_{s_t}^{i_t}...Y^{i_r}_{s_r}& if\ k=s_t,\ \ t\in\{1,...,r\}\\
0& if\ k\notin (s_1,...,s_r)
\end{array}
         \end{cases}.
\end{equation*}
If $k\notin(s_1,...,s_r),$ we get $F(w_i^k)=F(w)$ hence by Proposition \ref{aa}, 
\begin{equation*}
 \sum_{i=0}^dc_i^j(k)F(w_i^k)=\sum_{i=0}^dc_i^j(k)F(w)=0=D_k^jF(w)
\end{equation*}
Suppose now that $k\in (s_1,...,s_r)$ for example, let $k=s_r$ then
\begin{equation*}
 D_k^jF=\delta^{ji_r}\prod_{p=1}^{r-1}Y_{s_p}^{i_p}.
\end{equation*}
But,
\begin{eqnarray*}
 \sum_{i=0}^dc_i^j(k)F(w_i^k)&=&\sum_{i=0}^dc_i^j(k)\prod_{p=1}^rY_{s_p}^{i_p}(w_i^k)\\
&=&\left( \prod_{p=1}^{r-1}Y_{s_p}^{i_p}(w)\right)\sum_{i=0}^dc_i^j(k) Y_k^{i_r}(w_i^k)\\
&=&\left( \prod_{p=1}^{r-1}Y_{s_p}^{i_p}(w)\right)\sum_{i=0}^dc_i^j(k)v_i^{i_r}(k)\\
&=&\delta^{ji_r}\prod_{p=1}^{r-1}Y_{s_p}^{i_p}(w)\\
&=&D_k^jF(w).
\end{eqnarray*}
$\hfill\square$

An immediate consequence, is the following
\begin{cor}
 The gradient operator extends to any random variable $F:\Omega\longrightarrow\R.$
\end{cor}

We denote $Dom(D)$ the $L^2$-domain of $D$: $F\in Dom(D)$ if and only if
\begin{equation*}
 \E[\Arrowvert DF\Arrowvert^2_{l^2(\N)}]<\infty.
\end{equation*}

\section{Clark-Ocone formula}

In this section, we derive an explicit expression for the predictable representation of stochastic variables. 
The main tool is a discrete time analog of the well-known Clark-Ocone formula:  for any random variable $F$
\begin{equation*}
 F=\E[F]+\int \E[D_tF\,|\boF_t]dB_t.
\end{equation*}
When $d=1$, the discrete time analog of the Clark-Ocone formula for Bernoulli measures appears in \cite{P} 
(we refer also to \cite{L} for a discrete but finite Clark-Ocone formula). For general $d \geq 1$, one has
\begin{pro}
 For any $F\in\boS,$
\begin{eqnarray*}
 F&=&\E[F]+\sum_{k=0}^{\infty}<\E[D_kF\,|\boF_{k-1}],Y_k>\\
&=&\E[F]+\sum_{k=0}^{\infty}<D_k\E[F\,|\boF_{k}],Y_k>.
\end{eqnarray*}
\end{pro}
\pf By linearity, it suffices to show the result for $F=\boI^r(f_r),$ from the recurrence formula 
\begin{eqnarray*}
 F&=&r\sum_{i=1}^d\ \sum_{k=0}^{\infty}\boI^{r-1}\left( f_r^{i}(*,k)1_{\Delta_r}(*,k)1_{{\text{\textlbrackdbl}0,k-1\text{\textrbrackdbl}}^{r-1}}(*)\right) Y_{k}^{i}\\
&=&r\sum_{i=1}^d\ \sum_{k=0}^{\infty}\E\left[ \boI^{r-1}\left( f_r^{i}(*,k)1_{\Delta_r}(*,k)\right)\,|\boF_{k-1} \right] Y_{k}^{i}\\
&=&\sum_{i=1}^d\sum_{k=0}^{\infty}\E[D_k^iF\,|\boF_{k-1}]Y_k^i
\end{eqnarray*}
and since $\E[\boI^r(f_r)]=0,\ \forall r\geq1$ we get the first identity, while the second one holds from
\begin{equation*}
 \E[D_k^iF\,|\boF_{k-1}]=D_k^i\E[F\,|\boF_{k}].
\end{equation*}
$\hfill\square$

\begin{pro}
 The operator
\begin{equation*}
 \begin{array}{ccl}
L^2(\Omega)&\longrightarrow&L^2(\Omega\times\N,\R^d)\\
F&\longmapsto&((\E[D_k^1F\,|\boF_{k-1}],...,\E[D_k^dF\,|\boF_{k-1}]))_{k\in\N}
  \end{array}
\end{equation*}
is bounded with norm equal to one, hence the Clark-Ocone formula extends to any $F\in L^2(\Omega).$
\end{pro}
\pf From the Clark-Ocone formula and using the isometry formula, we have for $F\in\boS$;
\begin{eqnarray*}
 \Arrowvert\E[D_.F\,|\boF_{.-1}]\Arrowvert^2_{L^2(\Omega\times\N)}&=&\sum_{j=1}^d\Arrowvert\E[D_.^jF\,|\boF_{.-1}]\Arrowvert^2_{L^2(\Omega\times\N)}\\
&=&\sum_{j=1}^d\E\left[ \sum_{k=0}^{\infty}\left( \E[D_k^jF\,|\boF_{k-1}]\right)^2 \right] \\
&=&\E[(F-\E[F])^2]\\
&\leq&\E[\rvert F\rvert^2]\\
&=&\Arrowvert F\Arrowvert^2_{L^2(\Omega)}.
\end{eqnarray*}
$\hfill\square$\\
Consequently we state a Poincaré inequality.
\begin{cor}
 For any $F\in L^2(\Omega),$
\begin{equation*}
 var(F)\leq\Arrowvert DF\Arrowvert^2_{L^2(\Omega\times\N,\R^d)}.
\end{equation*}
\end{cor}
\pf We have
\begin{eqnarray*}
 var(F)&=&\E[\arrowvert F-\E[F]\arrowvert^2]\\
&=&\sum_{j=1}^d\E\left[\left( \sum_{k=0}^{\infty} \E[D_k^jF\,|\boF_{k-1}]\right)^2 \right] \\
&=&\sum_{j=1}^d\E\left[ \sum_{k=0}^{\infty}\left( \E[D_k^jF\,|\boF_{k-1}]\right)^2 \right] \\
&\leq&\sum_{j=1}^d\E\left[ \sum_{k=0}^{\infty} \E[\arrowvert D_k^jF\arrowvert^2\,|\boF_{k-1}] \right] \\
&=&\sum_{j=1}^d\E\left[ \sum_{k=0}^{\infty} \arrowvert D_k^jF\arrowvert^2 \right] \\
&=&\E\left[ \sum_{k=0}^{\infty} \Arrowvert D_kF\Arrowvert^2 \right]. 
\end{eqnarray*}
$\hfill\square$\\
Another variant of the Clark-Ocone formula is stated as
\begin{pro}\label{bb}
 For $n\in\N$ and $F\in L^2(\Omega)$. We have
\begin{equation*}
 F=\E[F|\,\boF_n]+\sum_{k=n+1}^{\infty}<\E[D_kF\,|\boF_{k-1}],Y_k>
\end{equation*}
and
\begin{equation*}
 \E[F^2]=\E[(\E[F\,|\boF_n])^2]+\E\left[\sum_{k=n+1}^{\infty}\|\E[D_kF\,|\boF_{k-1}]\|^2\right].
\end{equation*}
\end{pro}
\pf Since $\E[F\,|\boF_n]\in L^2(\Omega,\boF_n)$, the Clark-Ocone formula gives
\begin{eqnarray*}
 \E[F\,|\boF_n]&=&\E[F]+\sum_{i=0}^d\sum_{k=0}^n\E[D_k^j\E[F\,|\boF_n]\,|\boF_{k-1}]Y_k^j\\
&=&\E[F]+\sum_{i=1}^d\sum_{k=0}^nD_k^j\E[\E[F\,|\boF_n]\,|\boF_{k}]Y_k^j\\
&=&\E[F]+\sum_{i=1}^d\sum_{k=0}^nD_k^j\E[F\,|\boF_{k}]Y_k^j\\
&=&\E[F]+\sum_{i=1}^d\sum_{k=0}^n\E[D_k^jF\,|\boF_{k-1}]Y_k^j
\end{eqnarray*}
which proves the first identity.
The second identity is a consequence of the first one together with the isometry property of $\boI.$ 
$\hfill\square$

As an application of the Clark-Ocone formula, we obtain the following PRP
for discrete-time martingales.
\begin{pro}
 Let $(M_n)_{n\in\N}$ be a d-dimensional martingale in $L^2(\Omega)$ with respect to $(\boF_n)_{n\in\N}.$
There exists a predictable d-dimensional process $(\gamma_k)_{k\in\N}$ such that $\forall n\in\N,$
\begin{equation*}
 M_n^i=M_{-1}^i+\sum_{k=0}^n<\gamma_k^i,Y_k>,\qquad i\in\{1,...,d\}.
\end{equation*}
\end{pro}
\pf Let $k\geq 1.$ The Corollary \ref{bb} shows that $\forall i\in\{1,...,d\}$
\begin{eqnarray*}
 M_k^i&=&\E[M_k^i\,|\boF_{k-1}]+<\E[D_kM_k^i\,|\boF_{k-1}],Y_k>\\
&=&M_{k-1}^i+<\E[D_kM_k^i\,|\boF_{k-1}],Y_k>
\end{eqnarray*}
and one concludes by letting 
\begin{equation*}
 \gamma_k^i:=\E[D_kM_k^i\,|\boF_{k-1}],\quad k\geq0.
\end{equation*}
$\hfill\square$

\section{Divergence operator}

Let $\boU$ the subspace of $L^2(\Omega\times\N,\R^d)$ defined by
\begin{eqnarray*}
 \boU=\left\lbrace \sum_{r=0}^n\left( \boI^r\left(f_{r+1}^1(*,.)\right),...,\boI^r\left(f_{r+1}^d(*,.)\right)\right) ,\right. 
\\ \left. f_{r+1}\in L^2(\Delta_{r+1})\ \text{symmetric\  on\ its\  first\  $r$\ variables}\right\rbrace 
\end{eqnarray*}
\begin{Def}
 The divergence operator is the linear mapping $\delta:\boU\longrightarrow L^2(\Omega)$ defined by
\begin{equation*}
 \delta(X)=\delta\left( \boI^r\left(f_{r+1}^1(*,.)\right),...,\boI^r\left(f_{r+1}^d(*,.)\right)\right) =\boI^{r+1}\left( \widetilde{f}_{r+1}\right)
\end{equation*}
for $(X_k)_k$ of the form
\begin{equation*}
 X_k=\left( \boI^r\left(f_{r+1}^1(*,k)\right),...,\boI^r\left(f_{r+1}^d(*,k)\right)\right),\quad k\in \N.
\end{equation*}
\end{Def}
\begin{pro}
 The operator $\delta$ is the adjoint to $D$.
\end{pro}
\pf 
We consider $F=\boI^r(f_r)$ and $G=(G_k=(G_k^1,...,G_k^d))_{k\in\N}$ where $$G_k^i=\boI^s(g_{s+1}^i(*,k)), \forall 1\geq i\geq d.$$
We have
\begin{eqnarray*}
 &&\E\left[ <D_.F,G>_{l^2(\N)}\right] 
\\&&=\sum_{i=1}^d\E\left[ <D_.^iF,G^i>_{l^2(\N)}\right] \\
&&=\sum_{i=1}^dr\sum_{k=0}^{\infty}\E\left[ \boI^{r-1}(f_r^i(*,k)1_{\Delta_r}(*,k))\boI^s(g_{s+1}^i(*,k))\right] \\
&&=r!1_{\{r-1=s\}}\sum_{i=1}^d\sum_{
\begin{array}{l}
\scriptstyle i_1,...,i_{r-1}=1\\
\scriptstyle j_1,...,j_s=1
                    \end{array}
}^d\ <f_r^{i_1,...,i_{r-1},i},g_r^{j_1,...,j_{s},i}>1_{\{i_1=j_1,...i_{r-1}=j_s\}}\\
&&=\E\left[\boI^r(f_r)\boI^{s+1}(\widetilde{g}_{s+1})\right] \\
&&=\E\left[\delta(G)F\right] 
\end{eqnarray*}
$\hfill\square$\\
The following result shows that $\delta$ coincides with the stochastic integral operator $\boI$ on the square summable predictable processes.
\begin{pro}
 The operator $\delta$ can be extended to $L^2(\Omega\times\N,\R^d)$ with
\begin{equation}\label{Div}
 \delta(X)=\sum_{k=0}^{\infty}<X_k,Y_k> -\sum_{i=1}^d\sum_{k=0}^{\infty}<D_k(X_k^i),Y_k>Y_k^i,
\end{equation}
provided all series converge in $L^2(\Omega).$
\end{pro}
\pf Let $X_k=\left( \boI^r\left(f_{r+1}^1(*,k)\right),...,\boI^r\left(f_{r+1}^d(*,k)\right)\right),$ then one has
\begin{eqnarray*}
 &&\delta(X)\\ &&=\boI^{r+1}\left( \widetilde{f}_{r+1}\right)\\
&&=\sum_{k_1,...,k_{r+1}=1}^d\ \sum_{(i_1,...,i_{r+1})\in\Delta_{r+1}}\widetilde{f}_{r+1}^{k_1,...,k_{r+1}}(i_1,...,i_{r+1})Y_{i_1}^{k_1}...Y_{i_{r+1}}^{k_{r+1}}\\
&&=\sum_{k_1,...,k_{r},t=1}^d\ \sum_{k=0}^{\infty}\sum_{(i_1,...,i_{r})\in\Delta_{r}}
\widetilde{f}_{r+1}^{k_1,...,k_{r},t}(i_1,...,i_{r},k)Y_{i_1}^{k_1}...Y_{i_{r}}^{k_{r}}Y_k^t\\
&&-r\sum_{k_1,...,k_{r-1},s,t=1}^d\ \sum_{k=0}^{\infty}\sum_{(i_1,...,i_{r-1})\in\Delta_{r-1}}
\widetilde{f}_{r+1}^{k_1,...,k_{r-1},s,t}(i_1,...,i_{r-1},k,k)Y_{i_1}^{k_1}...Y_{i_{r-1}}^{k_{r-1}}Y_k^sY_k^t
\\&&=\sum_{k=0}^{\infty}\boI^{r}\left(f_{r+1}^t(*,k)\right)Y_k^t-r\sum_{s,t=1}^d\ \sum_{k=0}^{\infty}\boI^{r-1}\left(f_{r+1}^{s,t}(*,k,k)1_{\Delta_{r+1}}(*,k,k)\right)Y_k^sY_k^t
\\&&=\sum_{k=0}^{\infty}<X_k,Y_k> -\sum_{s,t=1}^d\ \sum_{k=0}^{\infty}D_k^s(X_k^t)Y_k^sY_k^t.
\end{eqnarray*}
$\hfill\square$\\
Observe that from  Proposition \ref{P}, the last term in the right-hand side of (\ref{Div}) vanish when $X$ is predictable.
As a result
\begin{cor} 
$\delta$ coincides with the stochastic integral operator $\mathcal{I}$ on the square summable predictable process.
\end{cor}

\section{Covariance identities}

The covariance $Cov(F,G)$ of $F,G\in L^2(\Omega)$ is defined as
\begin{equation*}
Cov(F,G)=\E[(F-\E[F])(G-\E[G])]=\E[FG]-\E[F]\E[G]. 
\end{equation*}

Let $(P_t)_{t\in\R_+}$ denote the Ornstein–Uhlenbeck semigroup, defined as
\begin{equation*}
 P_tF=\sum_{n=0}^{\infty}e^{-nt}\boI^n(f_n)
\end{equation*}
where $F=\sum_{n=0}^{\infty}\boI^n(f_n).$ We have
\begin{pro}\label{ci}
 For $F, G \in Dom(D),$
\begin{equation*}
 Cov(F,G)=\E\left[\sum_{k=0}^{\infty}<\E[D_kF\,|\boF_{k-1}],D_kG>\right]
\end{equation*}
and 
\begin{equation*}
 Cov(F,G)=\E\left[\sum_{k=0}^{\infty}\int_0^{\infty}e^{-t}<D_kF,P_tD_kG> dt\right].
\end{equation*}
\end{pro}
\pf The first identity is a consequence of the Clark-Ocone formula:
\begin{eqnarray*}
 Cov(F,G)&=&\E[(F-\E[F])(G-\E[G])]\\
&=&\E\left[\left(\sum_{k=0}^{\infty}<\E[D_kF\,|\boF_{k-1}],Y_k> \right) \left(\sum_{l=0}^{\infty}<\E[D_lG\,|\boF_{l-1}],Y_l> \right)  \right] \\
&=&\E\left[\left(\sum_{k=0}^{\infty}\sum_{i=1}^d\E[D_k^iF\,|\boF_{k-1}]Y_k^i \right) \left(\sum_{l=0}^{\infty}\sum_{i=1}^d\E[D_l^iG\,|\boF_{l-1}]Y_l^i\right) \right]\\
&=&\E\left[\sum_{k=0}^{\infty}\sum_{i=1}^d\E[D_k^iF\,|\boF_{k-1}]\E[D_k^iG\,|\boF_{k-1}]\right]\\
&=&\E\left[\sum_{k=0}^{\infty}\sum_{i=1}^d\E\left[\E[D_k^iF\,|\boF_{k-1}]D_k^iG\,|\boF_{k-1}\right]\right]\\
&=&\sum_{k=0}^{\infty}\E\left[\E\left[<\E[D_kF\,|\boF_{k-1}],D_kG>\,|\boF_{k-1}\right]\right]\\
&=&\E\left[\sum_{k=0}^{\infty}<\E[D_kF\,|\boF_{k-1}],D_kG>\right].
\end{eqnarray*}
By orthogonality of multiple
integrals of different orders and continuity of $P_t$ on $L^2(\Omega),$ it suffices to prove the second
identity for $F =\boI^n(f_n)$ and $G=\boI^n(g_n).$ We have
\begin{eqnarray*}
&&Cov(\boI^n(f_n),\boI^n(g_n))
\\&&=\E[\boI^n(f_n)\boI^n(g_n)]\\
&&= n!\sum_{k_1,...,k_n=1}^d<f_n^{k_1,...,k_n},g_n^{k_1,...,k_n}> \\
&&= n!n\int_0^{\infty}e^{-nt}dt\sum_{k_1,...,k_n=1}^d<f_n^{k_1,...,k_n},g_n^{k_1,...,k_n}> \\
&&= n\int_0^{\infty}e^{-nt}dt\sum_{k_1,...,k_n=1}^d\sum_{(i_1,...,i_n)\in\Delta_n}f_n^{k_1,...,k_n}(i_1,...,i_n)g_n^{k_1,...,k_n}(i_1,...,i_n)\\
&&= n^2\E\left[\int_0^{\infty}e^{-t} \sum_{j=1}^d \sum_{k=0}^{\infty}\boI^{n-1}\left( f_n^{j}(*,k)1_{\Delta_{n}}(*,k)\right)
e^{(n-1)t}\boI^{n-1}\left( g_n^{j}(*,k)1_{\Delta_{n-1}}(*,k)\right)dt\right] \\
&&= n^2\E\left[ \int_0^{\infty}e^{-t}\sum_{j=1}^d \sum_{k=0}^{\infty}\boI^{n-1}\left(f_n^{j}(*,k)1_{\Delta_{n}}(*,k)\right)
P_t\boI^{n-1}\left( g_n^{j}(*,k)1_{\Delta_{n-1}}(*,k)\right)dt\right]\\
&&= \E\left[ \int_0^{\infty}e^{-t}\sum_{j=1}^d \sum_{k=0}^{\infty}D_k^j\boI^{n}\left( f_n\right)
P_tD_k^j\boI^{n}\left( g_n\right)dt\right].
\end{eqnarray*}
$\hfill\square$

The next result shows that $(P_t)_{t\in\mathbb{R}_+}$ admits an integral representation by a probability kernel.

Let us define the probability kernel $Q_t(\widetilde{w},dw)$, for any $N\geq1$ and $t\in\mathbb{R}_+$, by
\begin{equation*}
 \mathbb{E}\left[ \left.\frac{dQ_t(\widetilde{w},.)}{d\mathbb{P}}  \right|\mathcal{F}_N\right] (w)=q_t^N(\widetilde{w},w)
\end{equation*}
where $q_t^N:\Omega\times\Omega\rightarrow \mathbb{R}_+$ is defined by
\begin{equation*}
 q_t^N(\widetilde{w},w)=\prod_{i=0}^N\left(1+e^{-t}<Y_i(w),Y_i(\widetilde{w})>\right),\quad w,\widetilde{w}\in\Omega.
\end{equation*}
\begin{pro}\label{pk}
 For any $F\in L^2(\Omega,\mathcal{F}_N)$ one has
\begin{equation*}
 P_tF(\widetilde{w})=\int_{\Omega}F(w)Q_t(\widetilde{w},dw),\quad \widetilde{w}\in\Omega.
\end{equation*}
\end{pro}
{\it Proof}: 
Recall that $L^2(\Omega,\mathcal{F}_N)$ has finite dimension $(d+1)^{N+1}$. More precisely, an orthonormal basis of $L^2(\Omega,\mathcal{F}_N)$ is given by
\begin{equation*}
 \left\lbrace Y^{s_1}_{k_1}...\,Y^{s_n}_{k_n}:\ 0\leq k_1<...<k_n\leq N,\ 1\leq s_1,...,s_n\leq d \right\rbrace.
\end{equation*}
Then it suffices to consider functionals of the form $F=Y_{k_1}^{s_1}...Y_{k_n}^{s_n}$.
Now observe that
\begin{eqnarray*}
 \mathbb{E}\left[ Y_k^j(.)\left(1+e^{-t}<Y_k(.),Y_k(w)>\right)\right] &=&\mathbb{E}\left[ Y_k^j(.)\left(1+e^{-t}\sum_{l=1}^dY_k^l(.)Y_k^l(w)\right)\right] 
\\&=&\sum_{i=0}^dc_i^j(k)(1+e^{-t}\sum_{l=1}^dv_i^l(k)Y_k^l(w))
\\&=&e^{-t}\sum_{i=0}^dc_i^j(k)\sum_{l=1}^dv_i^l(k)Y_k^l(w)
\\&=&e^{-t}\sum_{l=1}^d \delta^{jl}Y_k^l(w)
\\&=&e^{-t}Y_k^j(w).
\end{eqnarray*}
Then, by independence of the sequence $(X_k)_{k\geq0}$,
\begin{eqnarray*}
 \mathbb{E}\left[Y_{k_1}^{s_1}...Y_{k_n}^{s_n}q_t^N(\widetilde{w},.) \right]&=& \mathbb{E}\left[Y_{k_1}^{s_1}...Y_{k_n}^{s_n}\prod_{i=0}^N\left(1+e^{-t}<Y_{k_i}(\widetilde{w}),Y_{k_i}(.)>\right)\right]
\\&=& \prod_{i=0}^N\mathbb{E}\left[Y_{k_i}^{s_i}\left(1+e^{-t}<Y_{k_i}(\widetilde{w}),Y_{k_i}(.)>\right)\right]
\\&=&e^{-nt}Y_{k_1}^{s_1}(w)...Y_{k_n}^{s_n}(w)
\\&=&e^{-nt}\mathcal{I}^n\left( \widetilde{1}_{\{k_1,...,k_n\}}^{s_1,...,s_n}\right) (w)
\\&=&P_t(Y_{k_1}^{s_1}...Y_{k_n}^{s_n}).
\end{eqnarray*}
$\hfill\square$

\section{Deviation inequality}

In this section, we prove a deviation inequality for functionals of obtuse random walks, using the action of gradient operator and the covariance representations instead of the logarithm Sobolev inequality. We refer to \cite{P,HP} for other versions of this inequality in the one-dimensional case.
\begin{pro}\label{ED}
 Let $F:\Omega \rightarrow \mathbb{R}$ be a bounded random variable such that for any $k\in\mathbb{N}$ and $i,i'\in\{0,...,d\}$,
\begin{eqnarray*}
 |F(w_i^k)-F(w_{i'}^k)|\leq K,\quad |c_i^j(k)|\leq C,\ j\in\{1,...,d\}
\end{eqnarray*}
for some $K,C>0$ and $||DF||_{L^\infty(\Omega,l^1(\mathbb{N}))} <\infty$. Then
\begin{eqnarray*}
 \mathbb{P}(F-\mathbb{E}[F]\geq x)&\leq &\exp \left(-\frac{d\,C||DF||_{L^\infty(\Omega,l^1(\mathbb{N}))}}{K}g\left(\frac{x}{d\,C||DF||_{L^\infty(\Omega,l^1(\mathbb{N}))}} \right)  \right) 
\\&\leq& \exp \left(-\frac{x}{2K} \ln\left(1+\frac{x}{d\,C||DF||_{L^\infty(\Omega,l^1(\mathbb{N}))}} \right)  \right)
\end{eqnarray*}
with $g(u)=(1+u)\ln(1+u)-u,\ u\geq0$.
\end{pro}
{\it Proof}: For sake of simplicity, we assume that $\mathbb{E}[F]=0$. Then from the Chebychev inequality, one sees that
\begin{equation*}
 \mathbb{P}(F\geq x)\leq e^{-tx} \mathbb{E}\left[ e^{tF}\right] .
\end{equation*}
Next, letting $L(t)=\mathbb{E}\left[ e^{tF}\right]$, one has
\begin{eqnarray*}
 \ln(\mathbb{E}\left[ e^{tF}\right])&=&\int_{0}^t\frac{L'(s)}{L(s)}ds
\\&=&\int_{0}^t \frac{\mathbb{E}\left[ Fe^{sF}\right]}{\mathbb{E}\left[ e^{sF}\right]} ds.
\end{eqnarray*}
Now, we shall need the following result.
\begin{lem}\label{L1}
 For any random variable $F:\Omega \rightarrow \mathbb{R}$ and $s\geq0$, one has
\begin{equation*}
 e^{-sF}D_k^je^{sF}=\sum_{i=0, i\neq X_k}^dc_i^j(k)\left(e^{s(F(w_i^k)-F)}-1 \right) .
\end{equation*}

\end{lem}
{\it Proof}:
We have
\begin{eqnarray*}
 D_k^je^{F}&=&\sum_{i=0}^dc_i^j(k)e^{F(w_i^k)}
\\&=&\sum_{l=0}^d \textbf{1}_{\{X_k=l\}} \sum_{i=0,i\neq l}^d  c_i^j(k)e^{F(w_i^k)}+c_l^j(k)e^{F} 
\\&=&\sum_{l=0}^d \textbf{1}_{\{X_k=l\}} \sum_{i=0,i\neq l}^d  c_i^j(k)e^{F(w_i^k)}-\sum_{i=0,i\neq l}^d  c_i^j(k) e^{F} 
\\&=&\sum_{l=0}^d \textbf{1}_{\{X_k=l\}}e^F \sum_{i=0,i\neq l}^d c_i^j(k)\left( e^{F(w_i^k)-F}-1\right) .
\end{eqnarray*}
$\hfill\square$
\begin{rem}
 Note that the gradient operator does not satisfy a derivation rule for products. More precisely, for any $F,G:\Omega \rightarrow \mathbb{R}$, we have
\begin{equation*}
 D_k^j(FG)=FD_k^j(G)+GD_k^j(F)+\sum_{i=0, i\neq X_k}^dc_i^j(k)(F-F(w_i^k))(G-G(w_i^k)).
\end{equation*}

\end{rem}

Let $(P_t)_{t\in\mathbb{R}_+}$ be the Ornstein–Uhlenbeck semigroup, defined in the previous section.
From Proposition \ref{pk}, we obtain the following bound.
\begin{lem}\label{L3}
 For any $u\in L^2(\Omega\times\mathbb{N})$, one has
\begin{equation*}
 ||P_tu||_{L^\infty(\Omega,l^1(\mathbb{N}))}\leq ||u||_{L^\infty(\Omega,l^1(\mathbb{N}))}
\end{equation*}
\end{lem}
{\it Proof}: Using the representation formula of $P_t$ given in Proposition \ref{pk}, one has $\mathbb{P}(d\widetilde{w})$-a.s.
\begin{eqnarray*}
 ||P_tu||_{l^1(\mathbb{N})}(\widetilde{w})&=&\sum_{k=0}^\infty|P_tu_k(\widetilde{w})|
\\&=&\sum_{k=0}^\infty\left|\int_{\Omega}u_k(w)Q_t(\widetilde{w},dw)\right|
\\&\leq& \sum_{k=0}^\infty\int_{\Omega}|u_k(w)|Q_t(\widetilde{w},dw)
\\&=&\int_{\Omega}||u||_{l^1(\mathbb{N})}(w)Q_t(\widetilde{w},dw)
\\&\leq&||u||_{L^\infty(\Omega,l^1(\mathbb{N}))}.
\end{eqnarray*}
$\hfill\square$

Since the function $x\mapsto e^x-1$ is positive and increasing on $\mathbb{R}$, then by Lemma \ref{L1}, we obtain
\begin{eqnarray*}
 e^{-sF}D_k^je^{sF}&=&\sum_{i=0, i\neq X_k}^d c_i^j(k)\left(e^{s(F(w_i^k)-F)}-1\right)
\\&\leq& (e^{sK}-1)\sum_{i=0, i\neq X_k}^d c_i^j(k)
\\&=&- c_{X_k }^j(k) (e^{sK}-1)
 \\&\leq& C(e^{sK}-1).
\end{eqnarray*}
 
Going back to the proof of Proposition \ref{ED}, by Proposition \ref{ci} and Lemma \ref{L3} one has
\begin{eqnarray*}
 \mathbb{E}\left[ Fe^{sF}\right]&=& \text{Cov} (F,e^{sF})
\\&=&\mathbb{E}\left[\sum_{k=0}^{\infty}\int_0^{\infty}e^{-t}<D_ke^{sF},P_tD_kF> dt\right]
\\&\leq& d\,C(e^{sK}-1)\mathbb{E}\left[e^{sF}\int_0^{\infty}e^{-t}||P_tDF||_{(l^1(\mathbb{N}),\mathbb{R}^d)} dt\right]
\\&\leq& d\,C(e^{sK}-1)\mathbb{E}\left[e^{sF}\right]||DF||_{L^\infty(\Omega,l^1(\mathbb{N}))}\int_0^{\infty}e^{-t}dt
\\&\leq& d\,C(e^{sK}-1)\mathbb{E}\left[e^{sF}\right]||DF||_{L^\infty(\Omega,l^1(\mathbb{N}))}.
\end{eqnarray*}
Consequently we have
\begin{eqnarray*}
 \ln\left(\mathbb{E}\left[ e^{tF}\right]\right)&=&\int_{0}^t \frac{\mathbb{E}\left[ Fe^{sF}\right]}{\mathbb{E}\left[ e^{sF}\right]} ds
\\&\leq& d\,C||DF||_{L^\infty(\Omega,l^1(\mathbb{N}))}\int_0^t(e^{sK}-1)ds
\\&\leq& d\,C||DF||_{L^\infty(\Omega,l^1(\mathbb{N}))}(e^{tK}-tk-1).
\end{eqnarray*}
As a result, for any $x,t\geq0$,
\begin{equation*}
 \mathbb{P}(F\geq x)\leq \exp \left(d\,C||DF||_{L^\infty(\Omega,l^1(\mathbb{N}))}(e^{tK}-tk-1) -tx\right) .
\end{equation*}
The minimum in $t\geq0$ in the above expression is attained with
\begin{equation*}
 t=\frac{1}{K}\ln\left( 1+\frac{x}{d\,C||DF||_{L^\infty(\Omega,l^1(\mathbb{N}))}}\right) ,
\end{equation*}
whence
\begin{eqnarray*}
&& \mathbb{P}(F\geq x)
\\&&\leq \exp\left( -\frac{1}{K}((x+d\,C||DF||_{L^\infty(\Omega,l^1(\mathbb{N}))})\ln\left( 1+\frac{x}{d\,C||DF||_{L^\infty(\Omega,l^1(\mathbb{N}))}}\right) -x)\right) 
\\&&\leq\exp\left( -\frac{x}{2K}\ln\left( 1+\frac{x}{d\,C||DF||_{L^\infty(\Omega,l^1(\mathbb{N}))}}\right)\right) ,
\end{eqnarray*}
where we used the inequality $(1+u)\ln(1+u)-u\geq\frac{u}{2}\ln(1+u)$.
$\hfill\square$

\section{Complete markets in discrete time}

In this section we present a complete market model in discrete time 
as an application of the Clark-Ocone formula. 
The discrete-time and finite horizon models of financial markets can be described as follows.
We consider a probability space $(\Omega,\boF,(\boF_n)_{-1\leq n\leq N},\p)$ where $N$ is finite, $\boF_{-1}=\{\varnothing,\Omega\}$ and $\boF=\boF_N.$
Let 
\begin{equation*}
 \widetilde{S}=(B,S^1,...,S^d)=(B,S)\in\R^{d+1}
\end{equation*}
be $d+1$ assets such that $B=(B_n)_n$ is the price process of a risk-less asset (bond) where $B_n>0,\ \forall n\geq -1$
 and $S=(S^1,...,S^d)$ is the price process assets (stocks).
In order to distinguish the scalar product in $\R^{d+1}$, it will be convenient to use the notation $x.y=\sum_{i=1}^{d+1}x_iy_i$ for $x,y\in\R^{d+1}.$ 
We start by giving some classical definitions for financial market.
\begin{Def}
 A portfolio (or strategy) $\pi\in\R^{d+1}$ is a pair $(\beta,\gamma)$ , where $\beta=(\beta_n)$ and 
$\gamma_n=(\gamma_n^1,...,\gamma_n^d)$ are predictable processes such that
$\gamma_n^i$ denote the number of shares of the i$^{th}$ stock and $\beta_n$ the number of bonds that the seller of the option owns at time $n$.
\end{Def}
The corresponding value at time $n$ of the seller's portfolio is defined by
\begin{equation}\label{a}
 V_n^{\pi}=\pi_{n+1}.\widetilde{S}_n=\beta_{n+1} B_n+<\gamma_{n+1}, S_n>\qquad n\geq -1.
\end{equation}
 Note that in order to be consistent with the notation of the previous sections, we use the time scale $\N$, hence the index 0
is that of the first value of any stochastic process, while the index $-1$ corresponds to its deterministic initial value.
\begin{Def}
 A portfolio $\pi=(\beta,\gamma)$ is said to be self-financing if
\begin{equation}\label{b}
 B_{n}\Delta\beta_n+< S_{n},\Delta\gamma_n>=0,\qquad n\geq0,
\end{equation}
 
where $\Delta Y_n=Y_{n+1}-Y_{n}$ denotes the increment of the process $Y$ at time $n$.
\end{Def}
Hence the self-financing condition implies
\begin{equation*}
 V_n^{\pi}=\beta_{n} B_n+<\gamma_{n}, S_n>.
\end{equation*}
It's also convenient to use the discounted prices $\overline{S}=(\overline{S}^1,...,\overline{S}^d)$ defined by
\begin{equation*}
 \overline{S}_n= \frac{1}{B_n}S_n,\quad n\geq -1,
\end{equation*}
and the corresponding discounted value process of a strategy $\pi$ defined as
\begin{equation*}
 \overline{V}_n^{\pi}=\frac{1}{B_n}V_n^{\pi},\quad n\geq -1.
\end{equation*}
\begin{Def}
A model is said to be arbitrage-free if for every self-financing $\pi$ with $V_0^{\pi}=0$ and $V_N^{\pi}\geq0$ a.s.
then $V_N^{\pi}=0$ a.s.
\end{Def}
\begin{Def}
 An arbitrage-free market model is called complete if every contingent claim is attainable, i.e. every bounded 
$\boF$-measurable random variable $F$ can be hedged by a self-financing strategy.
\end{Def}
Note that in discrete-time setting, only a very limited class of models enjoys the completeness property. 
Let us recall the two basic theorems in asset pricing theory (see e.g. \cite{DMW,FS,JS} for proofs and more details).
\begin{Def}
 An equivalent martingale measure $\Q$ is a probability measure equivalent to $\p$  under which the d-dimensional discounted process $\overline{S}$ is a martingale.
\end{Def}

\begin{theo}\label{m}
 A market model is arbitrage-free if and only if there exists an equivalent martingale measure.
\end{theo}
\begin{theo}
 An arbitrage-free market model is complete if and only if there exists exactly one equivalent martingale measure. In this case, the number of atoms
in $(\Omega,\boF,\p)$ is bounded above by $(d+1)^N.$
\end{theo}
This theorem suggests that the price process of a complete market model,can be constructed from an obtuse random walk. 
Throughout the following, we are interested in market models given by
\begin{equation*}
\forall n\in\{0,...,N\},\quad
 \begin{cases}
 S_{n}=\prod_{k=0}^n(I+M_k^i)S_{-1}\quad \text{if}\ X_n=i,\ \ i\in\{0,...,l\}\\
B_{n}=\prod_{k=0}^n(1+r_k)
\end{cases}
\end{equation*}
with initial values $S_{-1}$ and $B_{-1}=1$, where $X_n$ is the coordinate maps i.e. 
\begin{equation*}
 \begin{array}{lccc}
   X_n:&\Omega&\longrightarrow&\{0,...,l\}\\
&w=(w_0,...,w_N)&\longmapsto&w_n
  \end{array}
\end{equation*}

and $(M_k^i)_k,\ (r_k)_k$ are deterministic sequences such that  $I+M_k^i$ is a matrix with
non-negative entries and $r_k>-1.$ 

\subsection{One-period model}

We start by discussing the notions of arbitrage-freeness and completeness of the market in one-period model
i.e. the assets are priced at the initial time $t=0$ and at the final time $t=1$. Let
\begin{equation*}
 \widetilde{\xi}=(1,\xi)=(1,\xi^1,...,\xi^d) \in\R^{d+1}_+
\end{equation*}
and
\begin{equation*}
 \widetilde{S}^i=(1,S^i)=(1+r,(I+M^i)\xi),\quad if\ X=i\in{0,...,l},
\end{equation*}
the respective asset price at time $t=0$ and $t=1$. Let us consider
\begin{equation*}
 \pi=(\pi^0,\gamma)=(\pi^0,\gamma^1,...,\gamma^d)\in\R^{d+1}
\end{equation*}
a portfolio at $t=0$. The price for buying $\pi$ equals $\pi.\overline{\xi}$ and
at time $t=1$, the portfolio takes the value
\begin{equation*}
\pi.\widetilde{S}^i=(1+r)\pi^0+<(I+M^i)\xi,\gamma>, 
\end{equation*}
depending on the scenario $X=i$.

With these notations, Theorem \ref{m} implies that the arbitrage-freeness is equivalent to the existence of
an equivalent martingale measure $\Q$ such that the probabilities $q_i=\Q(X=i)>0$ solve the linear equations
\begin{equation*}
 \begin{cases}
  q_0S^0+...+q_lS^l=r\xi\\
q_0+...+q_l=1. 
  \end{cases}
\end{equation*}
If a solution exists, it will be unique (i.e.  the arbitrage-free market model is complete) if and only if $l=d$ and
\begin{equation*}
 \binom{S^0}{1},...,\binom{S^d}{1}
\end{equation*}
 are linearly independent in $\R^{d+1}.$
\begin{rem}
 Note that if the arbitrage-free market model is complete, then the collection $\{S^0,...,S^d\}$ generates a convex set which contains the origin (c.f. \cite{FS}).
\end{rem}

\subsection{Multi-period model}

In the general case, the discounted prices $\overline{S}=(\overline{S}^1,...,\overline{S}^d)$ is given by
\begin{equation*}
 \overline{S}_n= \prod_{i=0}^n(1+r_i)^{-1}S_n,\quad n\geq -1,
\end{equation*}
and the corresponding discounted value process of a strategy $\pi$ defined by
 \begin{equation*}
  \overline{V}_n^{\pi}=\prod_{i=0}^n(1+r_i)^{-1}V_n^{\pi},\quad n\geq -1.
 \end{equation*}
The arbitrage-freeness is equivalent to the existence of an equivalent martingale measure $\Q$ such that 
the probabilities $q_n^i=\Q(X_n=i),\ n\in\N$ solve the linear equations
\begin{equation*}
 \begin{cases}
  q^0_{n+1}M_{n+1}^0S_n+...+q^l_{n+1}M_{n+1}^lS_n=r_{n+1}S_n\\
q^0_n+...+q^l_n=1. 
  \end{cases}
\end{equation*}

The arbitrage-free market model is complete if and only if $l=d$ and the matrix
\begin{equation*}
 \left( 
\begin{array}{ccc}
 M_{n+1}^0S_n& \cdots &M_{n+1}^dS_n  \\ 
1& \cdots &1
\end{array}
\right)\in\mathcal{M}_{d+1}(\R),
\end{equation*}
is invertible.
\begin{Ex}
 Consider
 \begin{equation*}
  M^i_k=diag(\lambda_k^{1,i},\ldots,\lambda_k^{d,i})\in \mathcal{M}_d(\R),\quad \text{for}\ k\in\{0,...,N\}\ \text{and}\ i\in \{0,...,d\}
 \end{equation*}
such that $\lambda_k^{j,i}>-1$ and the matrix
\begin{equation*}
 A_k=\left( 
\begin{array}{ccc}
 \lambda_k^{1,0}& \cdots &\lambda_k^{1,d}  \\ 
\vdots &  &\vdots  \\ 
 \lambda_k^{d,0}& \cdots& \lambda_k^{d,d} \\ 
1& \cdots &1
\end{array}
\right)\in\mathcal{M}_{d+1}(\R),
\end{equation*}
is invertible. 
\end{Ex}
Our goal is to provide a solution to the hedging problem : each $\boF$-measurable random variable $F$ can be hedged by a self-financing
strategy. In other terms, there is a self-financing $\pi$ such that
\begin{equation*}
 V_N^{\pi}=F\quad a.s.
\end{equation*}

In order to simplify the exposition, without losing much in generality, we assume that all $r_i$ are equal to $r$.
\begin{pro}
 The portfolio $\pi$ is self-financing if and only if,
\begin{equation}\label{c}
 V_n^{\pi}=V_{-1}^{\pi}+\sum_{k=0}^n\beta_{k}\Delta B_{k-1}+<\gamma_{k},\Delta S_{k-1}>.
\end{equation}
 
\end{pro}
\pf If $\pi$ is self-financing portfolio, it suffices to write
\begin{equation*}
 V_n^{\pi}=V_{-1}^{\pi}+\sum_{k=0}^nV_k^{\pi}-V_{k-1}^{\pi}.
\end{equation*}
Conversely, from (\ref{a}) we have
\begin{equation*}
 \Delta V_n^{\pi}= \beta_{n+1}B_n+<\gamma_{n+1},S_n>-\beta_{n}B_{n-1}-<\gamma_{n},S_{n-1}>.
\end{equation*}
But, the relation (\ref{c}) implies
 \begin{equation*}
  \Delta V_n^{\pi}= \beta_{n}\Delta B_{n-1}+<\gamma_{n},\Delta S_{n-1}>
 \end{equation*}
hence,
\begin{equation*}
 B_{n}\Delta\beta_n+< S_{n},\Delta\gamma_n>=0.
\end{equation*}
$\hfill\square$

\begin{pro}
 If the portfolio $\pi$ is self-financing, then
\begin{equation}\label{d}
 \Delta\overline{V}_n^{\pi}=<\gamma_{n+1},\Delta\overline{S}_n>,\quad n\geq-1.
\end{equation}
\end{pro}
\pf We write
\begin{eqnarray*}
 \Delta\overline{V}_n^{\pi}&=&(1+r)^{-n-1}V_{n}^{\pi}-(1+r)^{-n}V_{n-1}^{\pi}\\
&=&(1+r)^{-n-1}<\gamma_{n+1},S_n>-(1+r)^{-n}<\gamma_{n+1},S_{n-1}>.
\end{eqnarray*}
$\hfill\square$\\
The identity \ref{d} implies that
\begin{cor}
 
\begin{equation*}
 \overline{V}_n^{\pi}-\overline{V}_{n-1}^{\pi}=(1+r)^{-n-1}\sum_{j=1}^d(S_n^j-(1+r)S_{n-1}^j)\gamma_n^j.
\end{equation*}
\end{cor}
\begin{pro}
 Assume that the portfolio $\pi$ is self-financing. Then we have the decomposition
\begin{equation*}
 V_{n}^{\pi}=(1+r)^{n+1}V_{-1}^{\pi}+\sum_{k=0}^n\sum_{j=1}^d(1+r)^{n-k}(\lambda_k^{j,X_k}-r_k)\gamma_{k}^j S_{k-1}^j.
\end{equation*}
\end{pro}
\pf
We write
\begin{eqnarray*}
 V_{n}^{\pi}-V_{n-1}^{\pi}&=&\beta_{n}(B_{n}-B_{n-1})+<\gamma_{n},S_{n}-S_{n-1}>\\
&=&r\beta_{n}B_{n-1}+\sum_{j=1}^d\lambda_n^{j,X_n}\gamma_{n}^jS_{n-1}^j\\
&=&rV^{\pi}_{n-1}+\sum_{j=1}^d\lambda_n^{j,X_n}\gamma_{n}^jS_{n-1}^j.
\end{eqnarray*}
$\hfill\square$

The following proposition gives the unique equivalent martingale measure such that the market model is complete.

\begin{pro}
The process $(\overline{S}_n)_n$ is a d-dimensional martingale with respect to $(\boF_n)_{-1\leq n\leq N}$
under the probability $\Q$ given by
\begin{equation*}
 \left( \begin{array}{c}
       q_k^0\\
\vdots\\
q_k^d      
            \end{array}\right) = A_k^{-1}\left( \begin{array}{c}
       r\\
\vdots\\
r\\
1   
            \end{array}\right),\quad k\in\N.
\end{equation*}
Equivalently
\begin{equation*}
 \E_{\Q}[S_{n+1}\,|\boF_n]=(1+r)S_n,\qquad n\geq-1,
\end{equation*}
where $\E_{\Q}$ denotes the expectation under $\Q.$
\end{pro}

Recall that under this probability measure there is arbitrage freeness and the market is complete. By the predictable representation property
for discrete-time martingales we have
\begin{equation*}
 \overline{S}_n^i=\overline{S}_{-1}^i+\sum_{k=0}^n<\E_{\Q}[D_k\overline{S}_k^i\,|\boF_{k-1}],Y_k>.
\end{equation*}
So if $\pi$ is a self-financing portfolio, (\ref{d}) yields
\begin{eqnarray}
 \Delta\overline{V}_n^{\pi}&=&\sum_{i=1}^d\gamma_{n+1}^i<\E_{\Q}[D_n\overline{S}_n^i\,|\boF_{n-1}],Y_n>.
\end{eqnarray}
 
The next result provides a self-financing strategy $\pi$ as solution to the hedging problem.
\begin{pro}
 Let $F\in\boL^2(\Omega,\boF).$ For $ n\in\{0,...,N\}$, consider
\begin{equation}\label{f}
 \gamma_n^i=(1+r)^{n-N}\frac{1}{(S_n^i-(1+r)S_{n-1}^i)}\E_{\Q}[D_n^iF\,|\boF_{n-1}],\quad i\in\{1,...,d\},
\end{equation}
 and
\begin{equation}\label{g}
 \beta_n=(1+r)^{-N-1}\E_{\Q}[F\,|\boF_n]-(1+r)^{-n-1}<\gamma_n,S_n>.
\end{equation}
Then the portfolio $\pi=(\beta,\gamma)$ is self-financing and satisfies
\begin{equation*}
 V_n^{\pi}=(1+r)^{n-N}\E_{\Q}[F\,|\boF_n],\quad 0\leq n\leq N,
\end{equation*}
in particular $V_N^{\pi}=F$ a.s, hence $\pi$ is a hedging strategy leading to $F.$  
\end{pro}
\pf
Let $(\gamma_n)_{-1\leq n\leq N}$ be defined by (\ref{f}) with $\gamma_{-1}=0,$ and consider the process
$(\beta_n)_{-1\leq n\leq N}$ defined by
\begin{equation*}
 \begin{cases}
   \beta_{k+1}=\beta_k-(1+r)^{-k-1}<\Delta\gamma_k,S_k>,\quad k=-1,...,N-1\\
\beta_{-1}=(1+r)^{-N-1}\E_{\Q}[F]
  \end{cases},
\end{equation*}

such that  $\pi=(\beta_n,\gamma_n)_{-1\leq n\leq N}$ satisfies the self-financing condition. Let
\begin{equation*}
 \overline{V}_{-1}^{\pi}=(1+r)^{-N-1}\E_{\Q}[F],
\end{equation*}
hence, by the corollary (9.1), 
\begin{equation*}
 \overline{V}_n^{\pi}-\overline{V}_{n-1}^{\pi}=(1+r)^{-n-1}\sum_{j=1}^d(S_n^j-(1+r)S_{n-1}^j)\gamma_n^j.
\end{equation*}

On the other hand, from the Clark-Ocone formula, we have
\begin{eqnarray*}
&&(1+r)^{-N-1}\E_{\Q}[F\,|\boF_n]
\\&&=(1+r)^{-N-1}\E_{\Q}[F]+\sum_{k=0}^n\sum_{j=1}^d(1+r)^{-N-1}\E_{\Q}[D_k^jF\,|\boF_{k-1}]Y_k^j
 \\&&=(1+r)^{-N-1}\E_{\Q}[F]+\sum_{k=0}^n\sum_{j=1}^d(1+r)^{-k-1}(S_n^j-(1+r)S_{n-1}^j)\gamma_n^j.
\end{eqnarray*}
Hence
\begin{equation*}
 V_n^{\pi}=(1+r)^{n-N}\E_{\Q}[F\,|\boF_n].
\end{equation*}
In particular we have $V_N^{\pi}=F.$ Note that the relation (\ref{g}) follows from
\begin{equation*}
 V_n^{\pi}=\beta_{n} B_n+<\gamma_{n}, S_n>,\quad 0\leq n\leq N.
\end{equation*}
$\hfill\square$

We would like to thank an anonymous referee for helpful comments and suggestions.

\end{document}